\documentclass[a4paper,12pt]{article}
\usepackage{latexsym,amssymb}
\usepackage{amsmath}
\usepackage[mathscr]{eucal}
\usepackage{epic,eepic}
\usepackage{color}
\usepackage{enumerate}

\usepackage{amsthm}
\newtheorem{theorem}{Theorem}[section]

\newtheorem{corollary}[theorem]{Corollary}

\theoremstyle{definition}

\newtheorem{step}{Step}

\theoremstyle{remark}
\newtheorem{remark}[theorem]{Remark}

\makeatletter

\@addtoreset{equation}{section}
\makeatother

\newcommand{\ve}{\varepsilon}
\newcommand{\del}{\partial}
\newcommand{\lra}{\longrightarrow}
\newcommand{\e}{\mathrm{e}}
\renewcommand{\d}{\mathrm{d}}

\newcommand{\R}{\ensuremath{\mathbb{R}}}

\newcommand{\cL}{\mathcal{L}}
\newcommand{\sS}{\mathsf{S}}

\newcommand{\CAT}{\mathrm{CAT}}

\newcommand{\RCD}{\mathrm{RCD}}

\begin{document}

\title{How does the contraction property fail\\ for convex functions on normed spaces?}

\author{Shin-ichi OHTA\thanks{
Department of Mathematics, Osaka University, Osaka 560-0043, Japan
({\sf s.ohta@math.sci.osaka-u.ac.jp});
RIKEN Center for Advanced Intelligence Project (AIP),
1-4-1 Nihonbashi, Tokyo 103-0027, Japan}}

\date{\empty}
\maketitle

\begin{abstract}
On Euclidean and Hilbert spaces, Riemannian manifolds, and $\CAT(0)$-spaces,
gradient flows of convex functions are known to satisfy the contraction property,
which plays a fundamental role in optimization theory and possesses fruitful analytic and geometric applications.
On (non-inner product) normed spaces, however,
gradient flows of convex functions do not satisfy the contraction property.
We give a detailed proof of this characterization of inner products,
and discuss a possible form of a weaker contraction property on normed spaces.
%
\end{abstract}

\section{Introduction}

Analysis of gradient flows of convex functions is a fundamental subject in many branches
of mathematics, such as optimization theory, operator theory,
partial differential equations, geometric analysis, to name a few.
Beyond the classical setting of Euclidean and Hilbert spaces (see, e.g., \cite{Br}),
there has been a growing interest in gradient flows of convex functions on metric spaces.
(The convexity is then understood as the \emph{geodesic convexity} meaning that
the restriction of a function to every minimizing geodesic is convex.)
Such a generalization was intensively studied on \emph{$\CAT(0)$-spaces}
(nonpositively curved metric spaces in the sense of triangle comparison) in 1990s \cite{Jo1,Jobook,Ma};
we refer to \cite{AGSbook,Babook,Ba2} for further reading and more recent developments.
Besides $\CAT(0)$-spaces,
there are also a number of related studies on \emph{$\CAT(1)$-spaces} \cite{OP},
\emph{Alexandrov spaces of curvature bounded below}
as well as \emph{Wasserstein spaces} over them \cite{Ly,Ogra,Pe,Sa},
and metric measure spaces satisfying the \emph{Riemannian curvature-dimension condition}
$\RCD(K,\infty)$ \cite{St}.

On a metric space $(X,d)$, a gradient curve $\xi$ of a function $f$ is tending to the direction
in which $f$ decreases most efficiently (measured by the distance structure $d$),
which is a solution to the equation $\dot{\xi}(t) =\nabla[-f](\xi(t))$ in the smooth setting.
One of the most important properties of gradient curves $\xi,\zeta$ of a convex function is
the \emph{contraction} (or the \emph{non-expansion}) property:
\begin{equation}\label{eq:cont}
d\bigl( \xi(t),\zeta(t) \bigr) \le d\bigl( \xi(0),\zeta(0) \bigr) \quad \text{for all}\ t \ge 0.
\end{equation}
The contraction property plays a vital role from both theoretical and applied viewpoints.
It implies the continuous dependence (in particular, the uniqueness) of gradient curves on their initial points,
and gradient curves provide a contraction semigroup (\emph{gradient flow})
which can be analyzed by means of operator theory, partial differential equations, etc.

In the spaces mentioned in the first paragraph,
appropriately defined gradient curves of convex functions satisfy the contraction property \eqref{eq:cont}.
We remark and stress that those spaces are all ``Riemannian'' in the sense that
they exclude all non-Riemannian Finsler manifolds (and all non-inner product normed spaces).
In fact, on a normed space $(\R^n,\|\cdot\|)$,
we will show that the contraction property \eqref{eq:cont} holds for every convex function
if and only if the norm $\|\cdot\|$ comes from an inner product (Theorem~\ref{th:main}).
We remark that, in \cite[\S 4]{OSnc},
we have already seen an example of a convex function on a normed space
whose gradient flow does not have the contraction property.
We shall prove that, by testing only for the maximum of linear functions,
the contraction property fails for all non-inner product normed spaces.

The lack of the contraction property is one of the main difficulties of analysis of gradient flows
in ``non-Riemannian'' spaces.
There is no known quantitative estimate, and even the uniqueness of gradient curves is unclear
(see \cite[p.~4]{AGSbook} and \cite[Preface(b)]{Obook} for related remarks).
The main aim of this short article is to draw more attention to this fundamental open problem
and discuss a possible weaker contraction property available for normed spaces
and, more generally, Finsler manifolds (see Remark~\ref{rm:Q}).
Such a generalized contraction property will be obviously beneficial in optimization theory and Banach space theory,
as well as geometric analysis on Finsler manifolds, since heat flow on (measured) Finsler manifolds
can be regarded as gradient flow of the relative entropy on the \emph{$L^2$-Wasserstein space}
(we refer to \cite{OShf,OZ} for more details).
Moreover, solutions of some evolution equations can be regarded as gradient flow
in the $L^p$-Wasserstein space with $p \neq 2$,
which is non-Riemannian even when the underlying space is Riemannian
(see, for example, \cite{Ke} on $q$-heat flow).
We also refer to \cite{LOZ,Ohyp,Olln,OZ} and the references therein
for further related works concerning gradient flows in some classes of Finsler manifolds.

This article is organized as follows.
In Section~\ref{sc:pre}, we review necessary notions and recall the proof of the contraction property
in the Euclidean case.
Section~\ref{sc:main} is devoted to the proof of the failure of the contraction property for normed spaces
(Theorem~\ref{th:main}).
We also discuss a possible form of a weaker contraction property for normed spaces
inspired by Theorem~\ref{th:main}.

\section{Preliminaries}\label{sc:pre}

The canonical coordinates of $\R^n$ are denoted by $(x^i)_{i=1}^n$.
We will freely identify points and vectors in $\R^n$ in a canonical way.

\subsection{Gradient curves in normed spaces}\label{ssc:norm}

By a normed space $(\R^n,\|\cdot\|)$, we will mean a possibly asymmetric one; precisely,
\begin{enumerate}[(1)]
\item $\|x\| \ge 0$ for all $x \in \R^n$ and $\|x\|=0$ if and only if $x=0$;
\item $\|cx\|=c\|x\|$ for all $x \in \R^n$ and $c>0$;
\item $\|x+y\| \le \|x\| +\|y\|$ for all $x,y \in \R^n$.
\end{enumerate}
Then, the function $d(x,y):=\|y-x\|$ gives an asymmetric distance function on $\R^n$.
In this generality, a common notion of a \emph{gradient curve} for a locally Lipschitz function $f\colon \R^n \lra \R$
is a solution $\xi\colon [0,T) \lra \R^n$ to the \emph{energy dissipation identity}
\begin{equation}\label{eq:EDI}
f\bigl( \xi(t) \bigr) =f\bigl( \xi(s) \bigr) -\frac{1}{2} \int_s^t |\dot{\xi}|^2(r) \,\d r
 -\frac{1}{2} \int_s^t |\del f|^2\bigl( \xi(r) \bigr) \,\d r
\end{equation}
for all $0<s<t<T$, where
\[ |\dot{\xi}|(r) :=\limsup_{\rho \to r} \frac{\|\xi(\max\{r,\rho\})-\xi(\min\{r,\rho\})\|}{|r-\rho|} \]
is the \emph{$($forward$)$ metric derivative} of $\xi$ and
\[ |\del f|(x) :=\limsup_{y \to x} \frac{\max\{f(x)-f(y),0\}}{\|y-x\|} \]
is the \emph{local $($descending$)$ slope} of $f$.
Since ``$\ge$'' always holds in \eqref{eq:EDI} by the Cauchy--Schwarz and Young inequalities,
satisfying ``$\le$'' is the essential condition and implies ``$\dot{\xi}(r)=\nabla[-f](\xi(r))$'' at almost every $r$.
We refer to \cite{AG,AGSbook} for further reading,
and to \cite{OZ} and the references therein for the asymmetric case.

If $\|\cdot\|$ is strictly convex
(in the sense that $\|x+y\|=\|x\|+\|y\|$ holds only when $x=ay$ or $y=ax$ for some $a \ge 0$)
and $f$ is differentiable,
then we can introduce a gradient curve simply as a solution to
\[ \dot{\xi}(t) =\nabla[-f]\bigl( \xi(t) \bigr), \]
where $\nabla[-f](x) :=\cL(-\d f(x))$ is the \emph{gradient vector} defined via
the \emph{Legendre transform} $\cL\colon \R^n \lra \R^n$ associated with the norm $\|\cdot\|$ given by
\[ \cL(\alpha) := \biggl( \frac{1}{2} \frac{\del[\|\cdot\|_*^2]}{\del \alpha_i}(\alpha) \biggr)_{i=1}^n, \]
where $\|\cdot\|_*$ is the dual norm on $\R^n$, namely
\[ \|\alpha\|_* :=\sup\Biggl\{ \sum_{i=1}^n \alpha_i x^i \,\bigg|\, \|x\| \le 1 \Biggr\}. \]
(Note that $\|\cdot\|_*^2$ is differentiable since $\|\cdot\|$ is strictly convex.)
We remark that $\nabla[-f] \neq -\nabla f$ in general due to the asymmetry of $\|\cdot\|$.
In the case where $\|\cdot\|_*$ is twice differentiable on $\R^n \setminus \{0\}$ and $\d f(x) \neq 0$,
the gradient vector is written (with the help of Euler's homogeneous function theorem) as
\begin{equation}\label{eq:grad}
\nabla[-f](x) =\Biggl( -\sum_{j=1}^n \frac{1}{2}
 \frac{\del^2[\|\cdot\|_*^2]}{\del \alpha_i \del \alpha_j}\bigl( -\d f(x) \bigr) \frac{\del f}{\del x^j}(x) \Biggr)_{i=1}^n.
\end{equation}
We remark that $\|\cdot\|_*^2$ is twice differentiable at the origin $0$ only when $\|\cdot\|$ happens
to be induced from an inner product, thereby we excluded $0$
(see, for example, \cite[Proposition~1.7]{Obook}).

In \eqref{eq:grad},
\begin{equation}\label{eq:g_*}
g^*_{ij}(\alpha) :=\frac{1}{2} \frac{\del^2[\|\cdot\|_*^2]}{\del \alpha_i \del \alpha_j}(\alpha)
\end{equation}
(with $\alpha \neq 0$) is the inverse matrix of
\begin{equation}\label{eq:g_v}
g_{ij}\bigl( \cL(\alpha) \bigr) :=\frac{1}{2} \frac{\del^2[\|\cdot\|^2]}{\del x^i \del x^j}\bigl( \cL(\alpha) \bigr),
\end{equation}
provided that $(g_{ij}(\cL(\alpha)))_{i,j=1}^n$ is positive-definite.
By definition, for each $v \in \R^n \setminus \{0\}$, the symmetric matrix
$(g_{ij}(v))_{i,j=1}^n$ gives an inner product $g_v$ approximating $\|\cdot\|$ in the direction $v$
up to second order (see \cite[\S 3.1]{Obook} for a related account on Finsler manifolds).

\subsection{Contraction in Euclidean spaces revisited}\label{ssc:Eucl}

Before discussing the non-contraction in normed spaces,
let us recall how to show the contraction in the Euclidean setting.
We denote by $|\cdot|$ and $\langle \cdot,\cdot \rangle$
the standard Euclidean norm and inner product, respectively.

Let $f\colon \R^n \lra \R$ be a convex $C^1$-function
and $\xi,\zeta\colon [0,T) \lra \R^n$ be gradient curves of $f$.
For $t>0$, we have
\begin{align}
\frac{\d}{\d t} \Bigl[ |\zeta(t)-\xi(t)|^2 \Bigr]
&= 2\langle \zeta(t)-\xi(t),\dot{\zeta}(t)-\dot{\xi}(t) \rangle \label{eq:1vf}\\
&= -2\bigl\langle \zeta(t)-\xi(t),\nabla f \bigl( \zeta(t) \bigr) -\nabla f \bigl( \xi(t) \bigr) \bigr\rangle. \nonumber
\end{align}
Put $\gamma(s):=(1-s)\xi(t) +s\zeta(t)$ for $s \in [0,1]$.
Then the convexity of $f$ implies that
\[ f\bigl( \gamma(s) \bigr) \le (1-s)f\bigl( \xi(t) \bigr) +sf\bigl( \zeta(t) \bigr), \]
from which we find
\[ \frac{f(\gamma(s)) -f(\xi(t))}{s} \le \frac{f(\zeta(t)) -f(\gamma(s))}{1-s}. \]
By letting $s \to 0$ and $s \to 1$, since $\dot{\gamma}(0) =\dot{\gamma}(1) =\zeta(t)-\xi(t)$, we obtain
\begin{equation}\label{eq:mono}
\bigl\langle \nabla f\bigl( \xi(t) \bigr),\zeta(t)-\xi(t) \bigr\rangle
 \le f\bigl( \zeta(t) \bigr) -f\bigl( \xi(t) \bigr)
 \le \bigl\langle \nabla f\bigl( \zeta(t) \bigr),\zeta(t)-\xi(t) \bigr\rangle
\end{equation}
(that is, $\nabla f$ is monotone).
Hence, combining this with \eqref{eq:1vf}, we have the contraction
\[ \frac{\d}{\d t} \Bigl[ |\zeta(t)-\xi(t)|^2 \Bigr] \le 0. \]
Moreover, for a \emph{$K$-convex} function $f$ with $K \in \R$ in the sense that
\[ f\bigl( (1-s)x+sy \bigr) \le (1-s)f(x) +sf(y) -\frac{K}{2}(1-s)s|y-x|^2 \]
for all $x,y \in \R^n$ and $s \in (0,1)$, we have the \emph{$K$-monotonicity}
\[ \bigl\langle \nabla f\bigl( \xi(t) \bigr),\zeta(t)-\xi(t) \bigr\rangle
 -\bigl\langle \nabla f\bigl( \zeta(t) \bigr),\zeta(t)-\xi(t) \bigr\rangle
 \le -K|\zeta(t)-\xi(t)|^2 \]
as well as the \emph{$K$-contraction property}
\begin{equation}\label{eq:Kcont}
|\zeta(t)-\xi(t)| \le \e^{-Kt} |\zeta(0)-\xi(0)|.
\end{equation}

In the case of a normed space $(\R^n,\|\cdot\|)$, on the one hand, we have
\begin{align*}
\frac{\d}{\d t} \Bigl[ \|\zeta(t)-\xi(t)\|^2 \Bigr]
&= \sum_{j=1}^n \frac{\del[\|\cdot\|^2]}{\del x^j}\bigl( \zeta(t)-\xi(t) \bigr)
 \cdot \bigl( \dot{\zeta}_j(t)-\dot{\xi}_j(t) \bigr) \\
&= \sum_{i,j=1}^n \frac{\del^2[\|\cdot\|^2]}{\del x^i \del x^j}\bigl( \zeta(t)-\xi(t) \bigr)
 \cdot \bigl( \zeta_i(t)-\xi_i(t) \bigr) \bigl( \dot{\zeta}_j(t)-\dot{\xi}_j(t) \bigr)
\end{align*}
in place of \eqref{eq:1vf}, where we assumed that $\|\cdot\|$ is $C^2$ on $\R^n \setminus \{0\}$
as well as $\xi(t) \neq \zeta(t)$, and used the homogeneous function theorem in the second line.
Substituting $\dot{\xi}(t) =\nabla[-f](\xi(t))$ yields that,
with the help of \eqref{eq:grad}, \eqref{eq:g_*} and \eqref{eq:g_v},
\begin{align}
\frac{\d}{\d t} \Bigl[ \|\zeta(t)-\xi(t)\|^2 \Bigr] 
&= 2\sum_{i,j,k=1}^n g_{ij} \bigl( \zeta(t)-\xi(t) \bigr)
 \cdot \bigl( \zeta_i(t)-\xi_i(t) \bigr) \nonumber\\
&\quad \times
 \biggl( g^*_{jk} \bigl( -\d f \bigl( \xi(t) \bigr) \bigr) \frac{\del f}{\del x^k}\bigl( \xi(t) \bigr)
 -g^*_{jk} \bigl( -\d f \bigl( \zeta(t) \bigr) \bigr) \frac{\del f}{\del x^k}\bigl( \zeta(t) \bigr) \biggr).
 \label{eq:norm-1}
\end{align}
On the other hand, for $\gamma(s) =(1-s)\xi(t) +s\zeta(t)$ as above, we observe
\[ \lim_{s \to 0} \frac{f(\gamma(s)) -f(\xi(t))}{s}
 =\sum_{i=1}^n \frac{\del f}{\del x^i}\bigl( \xi(t) \bigr) \cdot \bigl( \zeta_i(t)-\xi_i(t) \bigr), \]
thereby the convexity of $f$ implies
\begin{equation}\label{eq:norm-2}
\sum_{i=1}^n \biggl( \frac{\del f}{\del x^i}\bigl( \xi(t) \bigr) -\frac{\del f}{\del x^i}\bigl( \zeta(t) \bigr) \biggr)
 \cdot \bigl( \zeta_i(t)-\xi_i(t) \bigr) \le 0
\end{equation}
in place of \eqref{eq:mono}.
Therefore, if
\begin{equation}\label{eq:norm-3}
\Bigl( g^*_{ij} \bigl( -\d f \bigl( \xi(t) \bigr) \bigr) \Bigr)
 =\Bigl( g^*_{ij} \bigl( -\d f \bigl( \zeta(t) \bigr) \bigr) \Bigr)
 =\Bigl( g_{ij} \bigl( \zeta(t)-\xi(t) \bigr) \Bigr)^{-1}
\end{equation}
holds as $n \times n$ matrices,
then we deduce from \eqref{eq:norm-1} and \eqref{eq:norm-2} the contraction
\[ \frac{\d}{\d t} \Bigl[ \|\zeta(t)-\xi(t)\|^2 \Bigr] \le 0. \]

When the norm $\|\cdot\|$ comes from an inner product,
$g_v$ (defined by $g_{ij}(v)$) coincides with the original inner product
for any $v \neq 0$ and \eqref{eq:norm-3} always holds.
On a normed space, however, $g_v$ depends on the direction $v$ and \eqref{eq:norm-3} does not hold
unless $\zeta(t)-\xi(t) =a\nabla[-f](\xi(t)) =b\nabla[-f](\zeta(t))$ for some $a,b>0$.

In view of \eqref{eq:norm-3}, the difficulty in the case of normed spaces
stems from the difference between the directions in question;
$\zeta(t)-\xi(t)$ from the first variation of the distance
and $\nabla[-f](\xi(t))$ (or $\nabla[-f](\zeta(t))$) from the gradient flow equation.
When \eqref{eq:norm-3} holds in some way, we can connect them
and apply the convexity \eqref{eq:norm-2} of $f$.
We refer to \cite{OP,OSnc} for related studies.
In \cite{OP} we introduced a condition called the \emph{commutativity} for metric spaces,
which corresponds to \eqref{eq:norm-3} and turned out a Riemannian condition.
In \cite{OSnc} we investigated a sufficient condition for the contraction,
called the \emph{skew convexity}.

\section{Non-contraction in normed spaces}\label{sc:main}

Let $(\R^n,\|\cdot\|)$ be a possibly asymmetric normed space in the sense of Subsection~\ref{ssc:norm}.
We will denote by $\sS :=\{x \in \R^n \mid \|x\|=1 \}$ the unit sphere of the norm.

We remark that it is natural to assume the strict convexity of $\|\cdot\|$.
In fact, without strict convexity, gradient curves are not unique
and there is no chance to have the contraction property
(see \cite[Example~4.23]{AG}, \cite[Remark~4.1]{OSnc}).

\begin{theorem}[Non-contraction in normed spaces]\label{th:main}
Suppose that a $($possibly asymmetric$)$ strictly convex normed space $(\R^n,\|\cdot\|)$
satisfies the \emph{contraction property} in the sense that
\begin{equation}\label{eq:contr}
\| \zeta(t)-\xi(t) \| \le \| \zeta(0)-\xi(0) \|
\end{equation}
holds for all $t \in [0,T)$ for every convex function $f\colon \R^n \lra \R$
and any gradient curves $\xi,\zeta\colon [0,T) \lra \R^n$ for $f$.
Then, the norm $\|\cdot\|$ is necessarily induced from an inner product.
\end{theorem}

We will see that, in fact, the contraction is tested only for the maximum of two linear functions in our proof.

\begin{proof}
We divide the proof into five steps.

\begin{step}\label{st:diff}
We first show that $\|\cdot\|$ is differentiable on $\R^n \setminus \{0\}$.
For a convex function, the directional differentiability in all directions implies the differentiability
(by \cite[Theorem~25.1]{Ro}).
Thus, assuming the contrary, let $\|\cdot\|$ be not directionally differentiable at $v \in \sS$.
Then, the closure of the set
\[ K := \{ x \in \R^n \mid v+\lambda (x-v) \in \sS\ \text{for some}\ \lambda >0 \} \]
is not a halfspace ($\del K$ is not a hyperplane) but a closed convex cone with origin at $v$.
Hence, we can modify the coordinates in such a way that
$v=(1,0,\ldots,0) \in \sS$ and $\|x\| >1$ holds for some $\delta >0$ and all
\[ x =(1+as,s,x^3,\ldots,x^n) \]
with $a \in (-\delta,\delta)$ and $(s,x^3,\ldots,x^n) \in \R^{n-1} \setminus \{0\}$
(to be precise, the strict inequality $\|x\|>1$ for all $(s,x^3,\ldots,x^n) \in \R^{n-1} \setminus \{0\}$
is guaranteed by the strict convexity).

\begin{figure}\label{fig:diff}
\centering\begin{picture}(240,160)

\put(120,0){\vector(0,1){160}}
\put(0,80){\vector(1,0){240}}
\put(150,25){$\sS$}

\qbezier(185,80)(170,110)(148,125)
\qbezier(148,125)(105,150)(80,125)
\qbezier(80,125)(60,105)(60,80)
\qbezier(60,80)(60,40)(80,25)
\qbezier(80,25)(105,5)(140,30)
\qbezier(140,30)(175,57)(185,80)

\thicklines
\put(165,160){\line(1,-4){40}}
\put(172,150){$x^1=1+ax^2$}

\put(185,90){\vector(1,0){30}}
\put(184,89){\rule{2pt}{2pt}}
\put(195,98){$\xi_+$}

\put(55,70){\vector(1,0){20}}
\put(54,69){\rule{2pt}{2pt}}
\put(64,57){$\xi_-$}

\end{picture}
\caption{Step~\ref{st:diff} for $n=2$}
\end{figure}
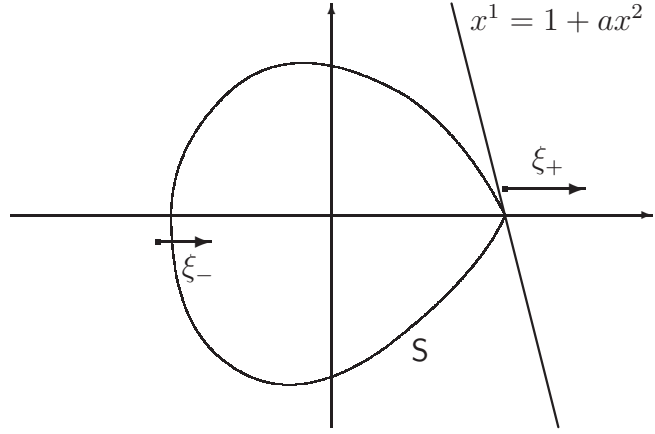

We consider the convex function
\[ f(x) =\max\{ -x^1+\ve x^2, c(-x^1-\ve x^2) \} \]
for $\ve \in (0,\delta)$ and $c \in (0,1)$.
By the choice of $\delta$, we have
\[ \nabla[-f](x) =\begin{cases}
 (1,0,\ldots,0) & \text{when}\ (1-c)x^1 < (1+c)\ve x^2, \\
 (c,0,\ldots,0) & \text{when}\ (1-c)x^1 > (1+c)\ve x^2. \end{cases} \]
Put
\[ b_c:= \frac{1-c}{1+c}\frac{2}{\ve} >0 \]
and observe
\[ (1+c)\ve b_c =2(1-c) >(1-c). \]
This implies that
\[ \nabla[-f]\bigl( (1,b_c,0,\ldots,0) \bigr) = (1,0,\ldots,0), \qquad
 \nabla[-f]\bigl( (-1,-b_c,0,\ldots,0) \bigr) = (c,0,\ldots,0). \]
Hence, the gradient curves emanating from $(1,b_c)$ and $(-1,-b_c)$ are given by
\begin{align*}
\xi_+(t) &= (1+t,b_c,0,\ldots,0) \quad \text{for}\ t \in [0,1), \\
\xi_-(t) &= (-1+ct,-b_c,0,\ldots,0) \quad \text{for}\ t \ge 0,
\end{align*}
respectively (see Figure~\ref{fig:diff}).
Thus, it follows from the assumed contraction property \eqref{eq:contr} that
\[\ \bigl\| \bigl( 2+(1-c)t,2b_c,0,\ldots,0 \bigr) \bigr\| \le \| (2,2b_c,0,\ldots,0) \| \]
for $t \in [0,1)$.
However, since $b_c \to 0$ as $c \to 1$, we have
\[ \bigl\| \bigl( 2+(1-c)t,2b_c,0,\ldots,0 \bigr) \bigr\| > \| (2,2b_c,0,\ldots,0) \| \]
when $c$ is close to $1$, a contradiction.
\end{step}

From here on, we assume that $\|\cdot\|$ is differentiable on $\R^n \setminus \{0\}$.
Note that $\|\cdot\|$ is then $C^1$ on $\R^n \setminus \{0\}$ (see \cite[Corollary~25.5.1]{Ro}).

\begin{step}\label{st:symm}

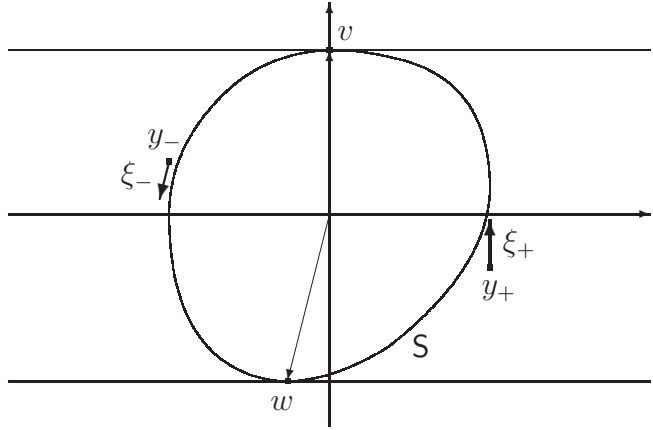
\begin{figure}\label{fig:symm}
\centering\begin{picture}(240,160)

\put(120,0){\vector(0,1){160}}
\put(0,80){\vector(1,0){240}}
\put(151,26){$\sS$}
\put(0,142){\line(1,0){240}}
\put(0,17){\line(1,0){240}}
\put(120,80){\vector(0,1){62}}
\put(120,80){\vector(-1,-4){15.6}}

\qbezier(180,90)(180,129)(148,138)
\qbezier(148,138)(105,150)(80,125)
\qbezier(80,125)(60,105)(60,80)
\qbezier(60,80)(60,40)(80,25)
\qbezier(80,25)(105,7)(142,30)
\qbezier(142,30)(180,60)(180,90)

\thicklines
\put(119,141){\rule{2pt}{2pt}}
\put(123,146){$v$}
\put(103.4,16){\rule{2pt}{2pt}}
\put(98,6){$w$}

\put(180,60){\vector(0,1){18}}
\put(179,59){\rule{2pt}{2pt}}
\put(177,50){$y_+$}
\put(185,66){$\xi_+$}
\put(60,100){\vector(-1,-4){3.6}}
\put(59,99){\rule{2pt}{2pt}}
\put(52,107){$y_-$}
\put(42,92){$\xi_-$}

\end{picture}
\caption{Step~\ref{st:symm} for $n=2$ ($\sS$ is $C^1$)}
\end{figure}

Next, to prove the symmetry of $\|\cdot\|$, we assume in contrary that the function
$\phi \colon \sS \lra (0,\infty)$, defined by the relation $-\phi(v) \cdot v \in \sS$,
is not constant (and hence its derivative is not identically zero).
Thus, by choosing appropriate coordinates,
the hyperplane $x^2=1$ is tangent to $\sS$ at $v:=(0,1,0,\ldots,0)$
and the hyperplane $x^2=-b$ is tangent to $\sS$ at $w:=(-a,-b,0,\ldots,0)$ for some $a,b>0$.
Now, given $c>0$, the function
\[ f(x)=\max\{ cb^{-1} x^2, -bx^2 \} \]
is convex and we have
\[ \nabla[-f](x) =
 \begin{cases} cw & \text{when}\ x^2>0, \\ bv & \text{when}\ x^2<0. \end{cases} \]
Hence, for $y_{\pm}:=(\pm 1,\mp s,0,\ldots,0)$ with small $s>0$,
\[ \xi_+(t) :=y_+ +tbv, \qquad \xi_-(t) :=y_- +tcw \]
are gradient curves emanating from $y_+, y_-$, respectively, for small $t \ge 0$ (see Figure~\ref{fig:symm}).
Thus, the hypothesized contraction property \eqref{eq:contr} implies that
\[ \frac{\d^+}{\d t} \Bigl[ \|\xi_+(t)-\xi_-(t)\| \Bigr]_{t=0}
 =\d\bigl[ \|\cdot\| \bigr]_{y_+ -y_-} (bv-cw) \le 0. \]
Letting $s \to 0$, since $\|\cdot\|$ is $C^1$ on $\R^n \setminus \{0\}$, we obtain
\[ \d\bigl[ \|\cdot\| \bigr]_{(2,0,\ldots,0)} (bv-cw) \le 0. \]
Moreover, by the contraction property between
\[ \zeta_+(t) :=(1,s,0,\ldots,0) +tcw, \qquad \zeta_-(t) :=(-1,-s,0,\ldots,0) +tbv, \]
we also find
\[ \d\bigl[ \|\cdot\| \bigr]_{(2,0,\ldots,0)} (cw-bv) \le 0. \]
Therefore,
\[ \d\bigl[ \|\cdot\| \bigr]_{(2,0,\ldots,0)} (bv-cw) =0 \]
holds for any $c>0$.
This means that $\sS$ is tangent to $(\lambda,0,\ldots,0) +\R \cdot (bv-cw)$
for $\lambda>0$ with $(\lambda,0,\ldots,0) \in \sS$ and arbitrary $c>0$.
However, this is possible only when $v$ and $w$ are linearly dependent, so that $a=0$.
This is a contradiction and completes the proof of the symmetry.
\end{step}

\begin{step}\label{st:square}
Given arbitrary $v \in \sS$, let us modify the coordinates in such a way that
$v=(1,0,\ldots,0)$ and the hyperplane $x^1=1$ is tangent to $\sS$ at $v$.
Then we consider the convex function
\[ f(x)=|x^1| \]
and shall show that lines parallel to $v$ and passing through $\sS \cap f^{-1}(0)$ are tangent to $\sS$
(this is indeed guessed from the argument in the previous step).

\begin{figure}\label{fig:square}
\centering\begin{picture}(240,160)

\put(120,0){\vector(0,1){160}}
\put(0,80){\vector(1,0){240}}
\put(60,0){\line(0,1){160}}
\put(180,0){\line(0,1){160}}
\put(138,28){$\sS$}

\qbezier(180,80)(180,129)(148,138)
\qbezier(148,138)(100,151)(74,119)
\qbezier(74,119)(60,100)(60,80)
\qbezier(60,80)(60,31)(92,22)
\qbezier(92,22)(140,9)(166,41)
\qbezier(166,41)(180,60)(180,80)

\thicklines
\put(179,79){\rule{2pt}{2pt}}
\put(184,71){$v$}

\put(119,140.7){\rule{2pt}{2pt}}
\put(109,131){$w$}
\put(179,94){\rule{2pt}{2pt}}
\put(184,95){$y_+$}
\put(59,64){\rule{2pt}{2pt}}
\put(44,65){$y_-$}

\put(180,95){\vector(-1,0){30}}
\put(160,103){$\xi_+$}
\put(60,65){\vector(1,0){30}}
\put(72,53){$\xi_-$}

\end{picture}
\caption{Step~\ref{st:square} for $n=2$ ($\sS$ is symmetric)}
\end{figure}
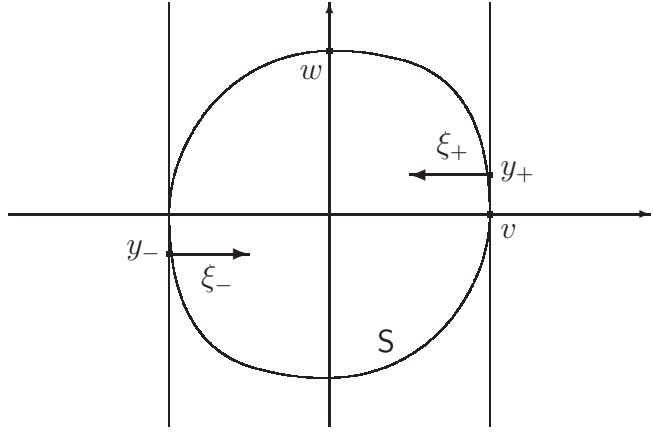

To this end, take $w \in \sS \cap f^{-1}(0)$ and put $y_+:=v+sw$ and $y_-:=-v-sw$ for $s>0$.
Note that $f(y_{\pm})=f(\pm v)=1$ and $\nabla[-f](y_{\pm}) =\nabla[-f](\pm v) =\mp v$,
and hence the gradient curves for $f$ emanating from $y_{\pm}$ are given by
\[ \xi_+(t) :=y_+ -tv, \qquad \xi_-(t) :=y_- +tv \]
for $t \in [0,1)$ (see Figure~\ref{fig:square}).
Then we observe from
\[ \| \xi_+(t)-\xi_-(t) \| =\| y_+-y_- -2tv \| =2\| v+sw-tv \| \]
that
\[ \frac{1}{2} \frac{\d}{\d t} \Bigl[ \| \xi_+(t)-\xi_-(t) \| \Bigr]_{t=0}
 = -\d\bigl[ \|\cdot\| \bigr]_{v+sw} (v)
 = -\d\bigl[ \|\cdot\| \bigr]_{(v/s)+w} (v)
 \to -\d\bigl[ \|\cdot\| \bigr]_w (v) \]
as $s \to \infty$,
where the latter equality is a consequence of the $1$-homogeneity of $\|\cdot\|$.
Hence, the contraction property \eqref{eq:contr} implies
\[ \d\bigl[ \|\cdot\| \bigr]_{w} (v) \ge 0. \]
Moreover, replacing $w$ with $-w$ gives the reverse inequality thanks to the symmetry of $\|\cdot\|$,
indeed,
\[ \d\bigl[ \|\cdot\| \bigr]_{-w} (v)
 =\d\bigl[ \|\cdot\| \bigr]_w (-v)
 =-\d\bigl[ \|\cdot\| \bigr]_w (v). \]
Therefore, we obtain
\[ \d\bigl[ \|\cdot\| \bigr]_w (v) =0. \]
This yields that the line $w+\R \cdot v$ is tangent to $\sS$ at $w$ as we claimed.
\end{step}

\begin{step}\label{st:ball}

\begin{figure}\label{fig:ball}
\centering\begin{picture}(240,160)

\put(120,0){\vector(0,1){160}}
\put(0,80){\vector(1,0){240}}
\put(140,29){$\sS$}

\qbezier(180,80)(180,129)(148,138)
\qbezier(148,138)(100,151)(74,119)
\qbezier(74,119)(60,100)(60,80)
\qbezier(60,80)(60,31)(92,22)
\qbezier(92,22)(140,9)(166,41)
\qbezier(166,41)(180,60)(180,80)

\thicklines
\put(60,0){\line(3,4){120}}
\put(13,0){$x^2=ax^1$}

\put(158.5,142){\rule{2pt}{2pt}}
\put(157,150){$y_+$}
\put(119,60){\rule{2pt}{2pt}}
\put(124,60){$y_-$}

\put(159.5,143){\vector(0,-1){9}}
\put(146,141){$\xi_+$}
\put(120,61){\vector(-1,0){12}}
\put(108,48){$\xi_-$}

\put(172,131){\rule{2pt}{2pt}}
\put(177,130){$y'_+$}
\put(107,79){\rule{2pt}{2pt}}
\put(103,86){$y'_-$}

\end{picture}
\caption{Step~\ref{st:ball} for $n=2$}
\end{figure}
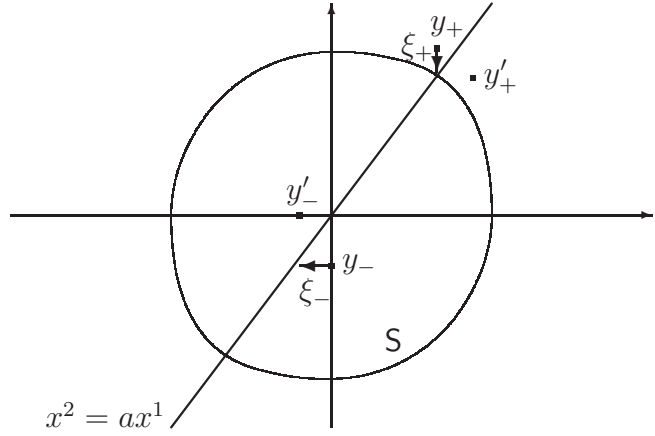

By virtue of the previous step, we can choose coordinates such that
the hyperplanes $x^1=1$ and $x^2=1$ are tangent to $\sS$
at $v=(1,0,\ldots,0)$ and $w=(0,1,0,\ldots,0)$, respectively.
Now, for $a>0$, we consider the convex function
\[ f(x) =\max\{ ax^1,x^2 \}. \]
Note that, by our choice of the coordinates,
\[ \nabla[-f](x) =\begin{cases} (-a,0,\ldots,0) & \text{when}\ ax_1>x_2, \\
 (0,-1,0,\ldots,0) & \text{when}\ ax_1<x_2. \end{cases} \]
We take $(r,ar,0,\ldots,0) \in \sS$ with $r=r(a)>0$ and put
$y_-:=(0,-as,0,\ldots,0)$, $y_+:=(r,ar+s,0,\ldots,0)$ for small $s>0$.
Observe that $\nabla[-f](y_-)=(-a,0,\ldots,0)$, $\nabla[-f](y_+)=(0,-1,0,\ldots,0)$, and
\begin{align*}
\xi_-(t) &:= y_- -(at,0) =(-at,-as,0,\ldots,0), \\
\xi_+(t) &:= y_+ -(0,t) =(r,ar+s-t,0,\ldots,0)
\end{align*}
are gradient curves for $f$ for $0\le t<\min\{s,s/a\}$ (see Figure~\ref{fig:ball}).
Then we have
\begin{align*}
\frac{\d}{\d t} \Bigl[ \|\xi_+(t) -\xi_-(t)\| \Bigr]_{t=0}
&= \d\bigl[ \|\cdot\| \bigr]_{(r,ar+(a+1)s)}(a,-1,0,\ldots,0) \\
&\to \d\bigl[ \|\cdot\| \bigr]_{(r,ar)}(a,-1,0,\ldots,0)
\end{align*}
as $s \to 0$.
Hence, by the contraction property \eqref{eq:contr}, we find
\[ \d\bigl[ \|\cdot\| \bigr]_{(r,ar)}(a,-1,0,\ldots,0) \le 0. \]
A similar discussion for $y'_- :=(-s,0)$ and $y'_+ :=(r+as,ar)$ yields that,
since $\nabla[-f](y'_-)=(0,-1,0,\ldots,0)$ and $\nabla[-f](y'_+)=(-a,0,\ldots,0)$,
\[ \d\bigl[ \|\cdot\| \bigr]_{(r,ar)}(-a,1,0,\ldots,0) \le 0. \]
Therefore, we obtain
\begin{equation}\label{eq:ball}
\d\bigl[ \|\cdot\| \bigr]_{(r,ar)}(a,-1,0,\ldots,0) =0,
\end{equation}
and the line $(r,ar,0,\ldots,0)+\R \cdot (a,-1,0,\ldots,0)$
is tangent to $\sS$ at $(r,ar,0,\ldots,0)$, for any $a>0$ and $r=r(a)>0$.
One can similarly show \eqref{eq:ball} for $a<0$ and $r=r(a)>0$ with $(r,ar,0,\ldots,0) \in \sS$.
Hence, the set consisting of $(b,c,0,\ldots,0) \in \sS$ coincides with the unit circle (i.e., $b^2 +c^2 =1$).
\end{step}

\begin{step}\label{st:last}
The outcome of the previous step is that
the parallelogram identity holds in the $2$-plane $P\langle v,w \rangle$ including $v$ and $w$.
By the way of taking $v,w$ in Step~\ref{st:square}, for any $\bar{v},\bar{w} \in \R^n$,
we can choose $v,w$ such that $\bar{v},\bar{w} \in P\langle v,w \rangle$.
Therefore, the parallelogram identity holds for any $\bar{v},\bar{w} \in \R^n$,
and hence $\|\cdot\|$ comes from an inner product.
This completes the proof.
\end{step}
\end{proof}

Comparing Theorem~\ref{th:main} with the scaling invariance of the convexity of functions on $\R^n$,
we can derive some immediate corollaries.
For this purpose, for a function $\Phi\colon (0,\infty) \lra (0,\infty)$,
we say that a strictly convex normed space $(\R^n,\|\cdot\|)$ satisfies the \emph{$\Phi$-contraction property} if
\[ \| \zeta(t)-\xi(t) \| \le \Phi(t) \| \zeta(0)-\xi(0) \| \]
holds for all $t \in (0,T)$ for every convex function $f\colon \R^n \lra \R$
and any gradient curves $\xi,\zeta\colon [0,T) \lra \R^n$ for $f$.

\begin{corollary}\label{cr:0}
Let $(\R^n,\|\cdot\|)$ be a $($possibly asymmetric$)$ strictly convex normed space,
which is not an inner product space.
If $(\R^n,\|\cdot\|)$ satisfies the $\Phi$-contraction property for some $\Phi\colon (0,\infty) \lra (0,\infty)$,
then we have
\[ \| \zeta(t)-\xi(t) \| \le \inf_{s>0} \Phi(s) \cdot \| \zeta(0)-\xi(0) \| \]
for all $t \in (0,T)$, for any convex function $f\colon \R^n \lra \R$
and gradient curves $\xi,\zeta\colon [0,T) \lra \R^n$ for $f$.
In particular, $(\R^n,\|\cdot\|)$ does not satisfy the $\Phi$-contraction property
for any $\Phi$ with $\inf_{s>0} \Phi(s) \le 1$.
\end{corollary}

Thus, the exponential contraction $\|\zeta(t) -\xi(t)\| \le \e^{-Kt} \|\zeta(0)-\xi(0)\|$
does not hold for any $K<0$ (compare this with \eqref{eq:Kcont}).

\begin{proof}
For any $\ve >0$, take $s>0$ such that $\Phi(s) \le \inf \Phi +\ve$.
Then, given a convex function $f\colon \R^n \lra \R$,
gradient curves $\xi,\zeta\colon [0,T) \lra \R^n$ for $f$, and $t \in (0,T)$,
we consider the function $cf$ with $c=t/s$.
Note that $cf$ is also convex and $\xi_c(t):=\xi(ct), \zeta_c(t):=\zeta(ct)$ are gradient curves for $cf$.
Hence, it follows from the $\Phi$-contraction property that
\begin{align*}
\|\zeta(t) -\xi(t)\|
&=\|\zeta_c(c^{-1}t) -\xi_c(c^{-1}t)\| \le \Phi(c^{-1}t) \|\zeta(0) -\xi(0)\| =\Phi(s) \|\zeta(0) -\xi(0)\| \\
&\le \bigl( \inf\Phi +\ve \bigr) \cdot \| \zeta(0)-\xi(0) \|.
\end{align*}
Since $\ve>0$ was arbitrary, we obtain the first assertion.
The latter assertion is then a consequence of Theorem~\ref{th:main}.
\end{proof}

\begin{remark}\label{rm:Q}
\begin{enumerate}[(a)]
\item
In view of Corollary~\ref{cr:0},
as a weaker contraction property on normed spaces,
what one could expect is a scaling invariant estimate of the form
\begin{equation}\label{eq:contC}
\|\zeta(t) -\xi(t)\| \le C \|\zeta(0) -\xi(0)\|
\end{equation}
with some $C>1$.
Then, due to the ``discontinuity'' at $t=0$,
it seems difficult to apply the method based on a differential inequality
(as in Subsection~\ref{ssc:Eucl} in the Euclidean case).

\item
Note that \eqref{eq:contC} suffices to show the uniqueness of gradient curves.
Thus, since the uniqueness fails in non-strictly convex normed spaces,
the constant $C$ would depend on the strength of the convexity of the norm $\|\cdot\|$,
such as the uniform convexity (see, for example, \cite[\S 1.2.1]{Obook}).

\item
Another direction is to consider $K$-convex functions for $K>0$.
Also in this case, we may need to have a closer look on the uniform convexity (and smoothness) of $\|\cdot\|$.
\end{enumerate}
\end{remark}
\medskip

\textit{Acknowledgements.}
This work was supported in part by JSPS Grant-in-Aid for Scientific Research (KAKENHI)
19H01786, 22H04942.



\end{document}